\documentclass[12pt]{article}

\usepackage{amsmath,amssymb}

\begin{document}

\pagestyle{myheadings} \markright{CLASS NUMBERS...}

\title{Asymptotics of class numbers}
\author{Anton Deitmar\footnote{Mathematisches Institut, Auf der Morgenstelle 10, 72076 T\"ubingen, Germany, \tt deitmar@uni-tuebingen.de} and Werner Hoffmann\footnote{Department
of Mathematical Sciences, South Road, Durham DH1 3LE, England;
\tt Werner.Hoffmann@durham.ac.uk}}

\date{}
\maketitle

$$ $$

\def \a{{{\mathfrak a}}}
\def \ad{\mathop{\rm ad}\nolimits}
\def \al{\alpha}
\def \ar{{\alpha_r}}
\def \A{{\mathbb A}}
\def \Ad{\mathop{\rm Ad}\nolimits}
\def \Aut{{\rm Aut}}
\def \b{{{\mathfrak b}}}
\def \bs{\backslash}
\def \B{{\cal B}}
\def \c{{\mathfrak c}}
\def \cent{\mathop{\rm cent}\nolimits}
\def \C{{\mathbb C}}
\def \CA{{\cal A}}
\def \CB{{\cal B}}
\def \CC{{\cal C}}
\def \CD{{\cal D}}
\def \CE{{\cal E}}
\def \CF{{\cal F}}
\def \CG{{\cal G}}
\def \CH{{\cal H}}
\def \CHC{{\cal HC}}
\def \CL{{\cal L}}
\def \CM{{\cal M}}
\def \CN{{\cal N}}
\def \CP{{\cal P}}
\def \CQ{{\cal Q}}
\def \CO{{\cal O}}
\def \CS{{\cal S}}
\def \CT{{\cal T}}
\def \CV{{\cal V}}
\def \d{{\mathfrak d}}
\def \det{\mathop{\rm det}\nolimits}
\def \df{\ \begin{array}{c} _{\rm def}\\ ^{\displaystyle =}\end{array}\ }
\def \diag{\mathop{\rm diag}\nolimits}
\def \dist{\mathop{\rm dist}\nolimits}
\def \End{\mathop{\rm End}\nolimits}
\def \eps{\varepsilon}
\def \eqn{\begin{eqnarray*}}
\def \endeqn{\end{eqnarray*}}
\def \F{{\mathbb F}}
\def \Fx{{\mathfrak x}}
\def \FX{{\mathfrak X}}
\def \fin{_{\mathrm{fin}}}
\def \g{{{\mathfrak g}}}
\def \ga{\gamma}
\def \Ga{\Gamma}
\def \Gal{{\rm Gal}}
\def \GL{\mathop{\rm GL}\nolimits}
\def \h{{{\mathfrak h}}}
\def \Hom{\mathop{\rm Hom}\nolimits}
\def \im{\mathop{\rm im}\nolimits}
\def \Im{\mathop{\rm Im}\nolimits}
\def \Ind{\mathop{\rm Ind}\nolimits}
\def \k{{{\mathfrak k}}}
\def \K{{\cal K}}
\def\Kfin{K_{\rm fin}}
\def \l{{\mathfrak l}}
\def \la{\lambda}
\def \li{{\rm li}}
\def \La{\Lambda}
\def \m{{{\mathfrak m}}}
\def \n{{{\mathfrak n}}}
\def \name{\bf}
\def \Mat{\mathop{\rm Mat}\nolimits}
\def \N{\mathbb N}
\def \o{{\mathfrak o}}
\def \ord{\mathop{\rm ord}\nolimits}
\def \O{{\cal O}}
\def \p{{{\mathfrak p}}}
\def \ph{\varphi}
\def \prf{\noindent{\bf Proof: }}
\def \Per{{\rm Per}}
\def \q{{\mathfrak q}}
\def \qed{\ifmmode\eqno $\square$\else\noproof\vskip 12pt plus 3pt minus 9pt \fi}
 \def\noproof{{\unskip\nobreak\hfill\penalty50\hskip2em\hbox{}%
     \nobreak\hfill $\square$\parfillskip=0pt%
     \finalhyphendemerits=0\par}}
\def \Q{\mathbb Q}
\def \res{\mathop{\rm res}\nolimits}
\def \R{{\mathbb R}}
\def \Re{\mathop{\rm Re}\nolimits}
\def \r{{\mathfrak r}}
\def \ra{\rightarrow}
\def \rank{\mathop{\rm rank}\nolimits}
\def \st{\mathop{\rm st}\nolimits}
\def \supp{\mathop{\rm supp}\nolimits}
\def \SL{\mathop{\rm SL}\nolimits}
\def \SO{\mathop{\rm SO}\nolimits}
\def \Spin{\mathop{\rm Spin}\nolimits}
\def \t{{{\mathfrak t}}}
\def \T{{\mathbb T}}
\def \tr{\mathop{\rm tr}\nolimits}
\def \vol{\mathop{\rm vol}\nolimits}
\def \z{\zeta}
\def \Z{\mathbb Z}
\def \={\ =\ }

\newcommand{\frack}[2]{\genfrac{}{}{0pt}{}{#1}{#2}}
\newcommand{\rez}[1]{\frac{1}{#1}}
\newcommand{\der}[1]{\frac{\partial}{\partial #1}}
\renewcommand{\binom}[2]{\left( \begin{array}{c}#1\\#2\end{array}\right)}
\newcommand{\norm}[1]{\left\|#1\right\|}
\renewcommand{\matrix}[4]{\left(\begin{array}{cc}#1 & #2 \\ #3 & #4 \end{array}\right)}
\renewcommand{\sp}[2]{\langle #1,#2\rangle}

\newtheorem{theorem}{Theorem}[section]
\newtheorem{conjecture}[theorem]{Conjecture}
\newtheorem{lemma}[theorem]{Lemma}
\newtheorem{corollary}[theorem]{Corollary}
\newtheorem{proposition}[theorem]{Proposition}

\section*{Introduction}

For an order $\CO$ in a number field let $h(\CO)$ denote its class number and
$R(\CO)$ its regulator.  Proving a conjecture of C.F. Gauss, C.L. Siegel showed in
\cite{siegel},
$$
\sum_{{d(\CO)\le x}}h(\CO) R(\CO)\ \sim \frac{\pi^2
}{36\zeta(3)} x^{\frac{3}{2}},
$$
where the sum ranges over the set of all 
real quadratic orders (i.e., orders in real quadratic fields) with discriminant
$d(\CO)$ bounded by
$x$.

For a long time it was
believed to be impossible to separate the class number and
the regulator. However, in 1981 P. Sarnak showed
\cite{sarnak}, using the trace formula, that
$$
\sum_{{R(\CO) \le x}}h(\CO) \ \sim \frac {e^{2x}}{2x},
$$
the sum ranging over all real quadratic orders with regulator bounded by $x$.
Sarnak established this result by identifying the regulators with
lengths of closed geodesics of the modular curve $H/\SL_2(\Z)$
(Theorem 3.1 there) and by using the prime geodesic
theorem for this Riemann surface. Actually, Sarnak proved not this result but the
analogue where $h(\CO)$ is replaced by the class number in the narrower sense and
$R(\CO)$ by a ``regulator in the narrower sense''. But in Sarnak's proof the
group $\SL_2(\Z)$ can be replaced by ${\rm PGL}_2(\Z)$ giving the above result. See
also \cite{Efrat, Venkov}.

Our goal is to generalize Sarnak's result to number fields of higher degree. Such a generalization has
resisted all efforts so far since the trace formula is not yet in a state that would make it useful in spectral
geometry. For instance, as yet a proof of the absolute convergence of the spectral side of the trace formula is
outstanding. However, recent partial results by W. M\"uller are
sufficient for the case treated in this paper.

We will now formulate the main theorem. Since there are several concepts of class
numbers, we have to make clear which one we use. Let $\CO$ be an order in a
number field
$F$. Let
$I(\CO)$ be the set of all finitely generated $\CO$-submodules of
$F$. According to the Jordan-Zassenhaus Theorem
\cite{reiner}, the set of isomorphism classes $[I(\CO)]$
of elements of $I(\CO)$ is finite. Let $h(\CO)$ be the
cardinality of the set $[I(\CO)]$, called the {\it class
number} of $\CO$. 

A cubic field $F$ (i.e., a number field of degree $3$ over the rationals) is either totally real or has two complex and one real embedding in
which case we call it a \emph{complex} cubic field.
Let $O$ be the set of isomorphism classes of orders in complex cubic
fields. The following is our main result.

\begin{theorem}\label{main}
As $x\ra\infty$ we have
$$
\sum_{\substack{\CO\in O\\R(\CO)\le x}} h(\CO)\ \sim\ \frac
{e^{3x}}{3x}.
$$
\end{theorem}

Our method is based on a new ``simple trace formula''. Such formulae have been used
in the past by various authors, for instance by Deligne-Kazhdan or Kottwitz. They
come about by plugging special test functions into Arthur's trace formula. The test
functions are chosen such that many terms in the trace formula vanish. 
The simple trace
formula of this paper will be such that the geometric side only consists of orbital
integrals of globally elliptic elements as opposed to locally elliptic
elements which is what the previous simple trace formulae reduced to.

In Section 1 the general form of the
simple trace formula is given, in which the test functions are characterized by
vanishing conditions. In Section 2, test functions are constructed explicitly by
twisting with virtual characters. This facilitates the computation of orbital
integrals and still leaves great freedom in the choice of test functions. The
convergence of the spectral side for noncompactly supported test functions is
discussed in Section 3. From this point on we restrict to the case $\SL_3$. In order
to separate orbital integrals of splitrank one, twisted resolvents are used. The
validity of the trace formula for these non-compactly supported functions is derived
via a Casimir functional calculus in Sections 4 to 7. The Prime
Geodesic Theorem, which is our main result in a different guise, is given in Section
8.

In the light of Sarnak's result and the result of the present paper, one is 
tempted to formulate a conjecture about the growth rate of class 
numbers in number fields of a given type. We must however warn that our 
method does not support any speculation of this kind. The reason is 
that we do not count orders, but rather their units which in the
setting of reductive groups come about as globally elliptic elements. Only in 
the case when the rank of the unit group is one, it is possible to draw
conclusions about class numbers from the distribution of units. Thus one is
limited to real  quadratic fields (Sarnak), complex cubic
(present paper) or purely  imaginary fields of degree 4. In the latter case a
new  difficulty emerges: since the degree of the field extension is not a
prime, ellitic elements can no longer be identified with order-units in number
fields any more, so the method gives a different asymptotic altogether.

{ \tableofcontents}

\section{A simple trace formula}\label{simple}

We derive a simple version of Arthur's trace formula by
inserting functions with certain restrictive properties which
guarantee the vanishing of the parabolic terms on the geometric side.
The trace formula for $SL(3,\Z)$ has also been studied in \cite{Wallace},
which unfortunately arrives at an incorrect formula due to a wrong
handling of the truncation.

Let $\CG$ be a linear algebraic $\Q$-group. If $E$ is a $\Q$-algebra,
any rational character  $\chi$ of $\CG$ defined over~$\Q$ defines
a homomorphism $\CG(E)\to\GL_1(E)$. If $E$ comes with an absolute
value $|\,.\,|$, we define $\CG(E)^1$ to be the subgroup of all
elements $g$ such that $|\chi(g)|=1$ for all
rational characters $\chi$ defined over~$\Q$. We will use this
notation in the cases when $E$ is $\R$ or the ring $\A$ of adeles
of~$\Q$. One should be aware that $\CG(\R)^1$ could also be
defined with respect to characters defined over the field~$\R$,
but this is not be the point of view in the present paper.

From now on we denote by $\CG$ a connected reductive linear algebraic
group over $\Q$. If $\CP$ is a parabolic $\Q$-subgroup of~$\CG$
with unipotent radical~$\CN$, we have a Levi decomposition
$\CP=\CL\CN$. Generally, we denote the group of real points of a
linear algebraic $\Q$-group by the corresponding roman letter, so that
$P=LN$. However, if $\CA$ is a maximal $\Q$-split torus of~$\CL$, we
denote by $A$ the connected component of the identity $\CA(\R)^0$.
One has decompositions $\CL(\A)=\CL(\A)^1A$, $L=MA$ (direct products)
and $P^1=MN$, where $M=L^1$.

Let $\A\fin$ denote the subring of finite adeles,
then we have direct product decompositions $\A=\A\fin \R$ and
$\CG(\A)=\CG(\A\fin)G$.
Fix a maximal compact subgroup $K$ of~$G$.

The geometric expansion of the trace formula
$$
J_{\mathrm{geom}}(f)=\sum_\o J_\o(f)
$$
for $f\in C_c^\infty(\CG(\A)^1)$ was introduced in \cite{Art-1}.
Here the sum runs over all classes
$\o$ in $\CG(\Q)$ with respect to the following equivalence relation:
Two elements are called equivalent if the semisimple components in their
Jordan decomposition are conjugate in~$\CG(\Q)$. Further,
$$
J_\o(f)=\int_{\CG(\Q)\bs\CG(\A)^1} k_\o(x)\,dx
$$
for certain functions $k_\o$ whose definition we will recall below. The sum and the integral converge if we replace $k_\o$ by its absolute value.

The integrand in the definition of $J_\o(f)$ is given as
$$
k_\o(x)=\sum_\CP (-1)^{\dim A_\CP/\CA_\CG}\sum_{\delta\in \CP(\Q)\bs\CG(\Q)}
K_{\CP,\o}(\delta x,\delta x)\hat{\tau}_\CP(H(\delta x)-T),
$$
where the sum runs over the standard parabolic $\Q$-subgroup~$\CP=\CL\CN$, for which we write $\CA_\CP=\CA_\CL$ and
$$
K_{\CP,\o}(x,y) \= \sum_{\ga\in\CL(\Q)\cap\o} \int_{\CN (\A )}f(x^{-1}\ga ny)\ dn.
$$
All we need to know about the factor $\hat\tau_\CP(H(\delta x)-T)$ at this point is that in the case
$\CP=\CG$ it is identically equal to~1.

We call a function on $\CG(\A)^1$ {\it parabolically regular
at the infinite place\/} if it is supported on $\Kfin \times G^1$ for some
compact open subgroup $K\fin$ of $\CG(\A\fin)$ and vanishes on all
$G$-conjugates of $\Kfin \times P^1$ for every parabolic
$\Q$-subgroup $\CP\ne\CG$.

An element $\gamma\in\CG(\Q)$ is called {\it $\Q$-elliptic\/} if it
is not contained in any parabolic $\Q$-subgroup other than $\CG$
itself. This notion is clearly invariant under conjugation, and we say
that a class $\o$ is $\Q$-elliptic if some (hence any) of its elements
is so. It is known that $\Q$-elliptic elements are semisimple, so
$\Q$-elliptic classes $\o$ are just conjugacy classes in~$\CG(\Q)$.

\begin{proposition} \label{simple-trace-formula}
If $f\in C_c^\infty(\CG(\A)^1)$ is parabolically regular at the infinite place,
then $J_\o(f)$ vanishes unless $\o$ is $\Q$-elliptic, in which case
$$
J_\o(f)=\int_{\CG(\Q)\bs\CG(\A)}K_{\CG,\o}(x,x)\,dx.
$$
\end{proposition}

In light of the above remarks the proposition is a consequence of the following lemma.

\begin{lemma}\label{1.2}
Suppose that $f$ satisfies the conditions of Proposition~\ref{simple-trace-formula}. Then $K_{\CP,\o}(x,x)=0$ for any $x\in\CG(\A)$ unless $\CP=\CG$ and
the class~$\o$ is $\Q$-elliptic.
\end{lemma}

\prf
First we show for any $\CP\ne\CG$ that $f(x^{-1}qx)=0$ for $q\in \CP(\A )^1$
and $x \in \CG(\A)$.
By the assumption on the support of $f$ we have only to consider
$q=q\fin q_\infty$ with $x^{-1}q\fin x \in \Kfin $, i.e.,
$q\fin \in x \Kfin  x^{-1} \cap\CP(\A )$, a compact subgroup of $\CP (\A )$.
Any continuous quasicharacter with values in $]0,\infty[$ will be
trivial on that subgroup, hence $q\fin \in \CP (\A )^1$.
Since $q$ was already in $\CP(\A )^1$, it follows that 
$q_\infty\in\CP(\A)^1\cap P =P^1$, and so $f(x^{-1}qx)=0$ due to the
assumption on $f$ applied to the parabolic $\Q$-subgroup $\CP$.

In particular, as $\CL(\Q)\CN(\A)\subset\CP(\A)^1$, we see that
$f(x^{-1}\ga nx)=0$ for all $\ga\in\CL(\Q)$ and $n\in \CN(\A )$, hence
$K_{\CP,\o}(x,x)=0$.

If $\o$ is not $\Q$-elliptic, then every $\ga\in\o$ is contained in some
parabolic $\Q$-subgroup $\CP\ne\CG$, and in view of $\CP(\Q)\subset\CP(\A)^1$
we have $f(x^{-1}\ga x)=0$. Thus $K_{\CG,\o}(x,x)=0$. Lemma \ref{1.2} and Proposition \ref{simple-trace-formula} follow.
\qed

We will now rewrite, in a non-adelic language, the geometric side of our simple trace formula in a special case. 
Let $\CG$ be a semisimple simply-connected linear algebraic $\Q$-group such
that $G=\CG(\R)$ has no compact factors.  Let  $\Gamma$ be a \emph{congruence
subgroup} of
$\CG(\Q)$, i.e., assume that there exists an open compact subgroup $K_\Ga$ of
$\CG(\A\fin)$ with
$\Ga= K_\Ga\cap\CG(\Q)$.  Let $f_\infty\in C^\infty_c(G)$ be a parabolically
regular function. For $y\in G$ the {\it orbital integral\/} is defined as
$$
\CO_y(f_\infty)\ :=\ \int_{G_y\bs G}f_\infty(x^{-1} yx)\ dx,
$$
where $G_y$ denotes the centralizer of $y$ in $G$. We  define $f=f\fin\otimes f_\infty$, where $f\fin$ is the characteristic function of $K_\Gamma$ divided by the volume of that group with respect to the Haar measure of~$\CG(\A\fin)$.

\begin{corollary}\label{geomtrace} Under the above conditions, we have
$$
J_{\mathrm{geom}}(f) \= \sum_{[\ga]}\ \vol (\Ga_\ga \bs G_\ga)\ \CO_\ga(f_\infty),
$$
where the sum on the right-hand side runs over the set of all conjugacy classes
$[\ga]$ in the group $\Ga$ which consist of $\Q$-elliptic elements.
\end{corollary}

\prf
Consider the formula for $J_\o(f)$ given in Proposition~\ref{simple-trace-formula} for a $\Q$-elliptic class $\o$. The integral can be taken over $\CG(\Q)\bs\CG(\A )/{K_\Ga}$, because
the integrand is right $K_\Gamma$-invariant.
Under our assumptions on~$\CG$, strong approximation \cite{Kn} holds,
i.e., the action of $G$ by right translation on that double quotient
is transitive and hence induces an isomorphism of $G$-spaces
$$
\Ga \bs G  \stackrel{\sim}{\longrightarrow}
\CG(\Q)\bs\CG(\A )/{K_\Ga}.
$$
This isomorphism identifies suitably normalized $G$-invariant measures on these two spaces  with each other, and we get
$$
J_\o(f)=\int_{\Gamma\bs G}\sum_{\ga\in\o}f\fin(\ga)f_\infty(x^{-1}\ga x)\,dx.
$$
The characteristic function $f\fin$ has the effect of restricting summation to $\o\cap\Gamma$.
% It is known that this set is a finite union of $\Gamma$-conjugacy classes.
Note that $J_\o(|f|)<\infty$, because $|f_\infty|$ can be bounded by a nonnegative function in $C_c^\infty(G)$.
Now the equality of $J_\o(f)$ with the partial sum over $[\ga]\in\o$
in the asserted formula follows by the familiar Fubini-type argument.
\qed

\section{Test functions}
In this section $G$ is a semisimple real Lie group with finite center and finitely
many connected components.

We would like to use resolvent kernels as test functions in the trace formula. The
convergence of the geometric side has been proved in~\cite{Art-1} for compactly
supported test functions only. Thus we are going to approximate resolvent kernels
by compactly supported functions with the aid of the functional calculus of the
Casimir operator. 

Let $\g_\R$ denote the Lie algebra of $G$ and $\g=\g_\R\otimes\C$ its
complexification. Let $B$ denote the Killing form. Let $\theta$ be the Cartan
involution fixing
$K$. The form $\sp XY=-B(\theta(X),Y)$ is positive definite on $\g_\R$ and induces a $G$-invariant Riemannian metric on $G/K$. Let
$\dist(x,y)$ denote the corresponding distance function and write $d(g)=\dist(gK,eK)$ for $g\in G$.

Let $U(\g)$ denote the universal enveloping algebra of $\g$. Every element $X$
of $U(\g)$ gives rise to a left-invariant differential operator written
$h\mapsto h * X$, and a right-invariant differential operator $h\mapsto X *
h$, $h\in C^\infty(G)$. Recall that for
$p>0$, the
$L^p$-Schwartz space
$\CC^p(G)$ is defined as the space of all $h\in C^\infty(G)$ such that, for
every $n>0$ and
$X$,
$Y\in U(\g)$, the seminorms
\[
|h|_{p,n,X,Y}=\sup_{g\in G}|X*h*Y(g)|\,\Xi(g)^{-2/p}(1+d(g))^n
\]
are finite. Here $\Xi$ is the basic spherical function, and it suffices for our present purposes to know that there exist $r_1>r_2>0$ such that $e^{-r_1d(g)}\le|\Xi(g)|\le e^{-r_2d(g)}$. If we complete the space $\CC^p(G)$ with respect to the seminorms involving only derivatives up to order~$N$, we obtain a space $\CC^p_N(G)$, whose topology can be given by a Banach norm.
For each $\tau\in\hat K$, there is a subspace $\CC^p_N(G,\tau)$ of functions $h$ satisfying $\chi_\tau*h*\chi_\tau=h$, where $\chi_\tau\in C(K)$ is the idempotent associated to~$\tau$.

We also need the space $\CH^r_N$ of even holomorphic functions $\phi$ on the strip $\{z\in\C\mid|\Im z|<r\}$ extending continuously to the boundary and such that the norm
\[
|\phi|_{r,N}=\sup_{|\Im z|\le r}|\phi(z)|(1+|\Re z|)^N
\]
is finite. Recall that a Schwartz function $\phi$ on $\R$ is called a \emph{Paley-Wiener function} if its Fourier transform
$\hat\phi(x)=\frac 1{2\pi}\int_\R\phi(y)\,e^{-ixy}\,dy$ has compact support.

For $\pi\in\hat G$ and an irreducible unitary representation
$(\tau,V_\tau)$ of $K$, let $P_{\pi,\tau}$ be the orthogonal projection
defined on the space of $\pi$ whose image is
the $\tau$-isotypical component.
Let $C$ be the Casimir operator of $G$.
For $\pi\in\hat G$ the Casimir $C$ acts on $\pi$ by a scalar $\pi(C)$.

\begin{proposition}\label{test}
Let $0<p\le1$, $N\in\N$, $b\in\R$ and $\tau\in\hat K$ be given. Then there exist $r>0$ and $N'\in\N$ such that for every $\phi\in\CH_{N'}^r$ there is a unique function $h_{\phi,\tau}\in \CC^p_N(G,\tau)$ satisfying
$$
\pi(h_{\phi,\tau})\=
\frac1{\dim\tau}\,\phi\left(\sqrt{-\pi(C)-b}\right) P_{\pi,\tau}
$$
for every $\pi\in\hat G$. The map  $\CH^r_{N'}\to\CC^p_N(G,\tau)$ so defined is continuous. If $\phi$ is a Paley-Wiener function, then $h_{\phi,\tau}$ is compactly supported.
\end{proposition}

(The factor $\frac 1{\dim\tau}$ is put here in order to give
$\pi(h_{\phi,\tau})$  a nice trace.)

\prf
The uniqueness of $h_{\phi,\tau}$ is clear from the Plancherel theorem. 
In the case of an even Paley-Wiener function~$\phi$, the existence of $h_{\phi,\tau}\in C_c^\infty(G)$ with the required properties (except for the bounds) has been proved in~\cite{limitmult}, Lemmas 2.9 and~2.11. 
For this one considers the $G$-homogeneous vector bundle  $E_\tau=G\times V_\tau/K$ and identifies the space of smooth sections with $\left(C^\infty(G)\otimes V_\tau\right)^K$. 
Then the operator
$D_\tau$ induced by $-C-b$ is a \emph{generalized Laplacian} in the sense of \cite{BGV}. The operator
$\phi\left(\sqrt{D_\tau}\right)$ defined by functional calculus is
an operator with smooth kernel
$\langle x | \phi\left(\sqrt{D_\tau}\right) | y\rangle$, and by the theory of
hyperbolic equations (\cite{Tayl}, ch.~IV) it follows that
$\phi\left(\sqrt{D_\tau}\right)$ has finite propagation speed (compare
\cite{CGT}). Identifying the sections of $E_\tau$ with $K$-invariant
functions as above, it follows that the $G$-equivariant operator
$\frac 1{\dim\tau}\,\phi\left(\sqrt{D_\tau}\right)$ acts as a convolution operator by a function $h_{\phi,\tau}$.

Now observe that by the estimates in \cite{CGT} there exists a constant $c>0$ such that, for every $\phi$,
$$
\bigl|\bigl\langle x\bigl|\phi\bigl(\sqrt{D_\tau}\bigr)\bigr|y\bigr\rangle\bigr|\ \le\
c|\hat\phi(\dist(x,y))|.
$$
This means that there is a constant $c>0$ such that for every $\phi$ we have
$$
|h_{\phi,\tau}(g)|\ \le\ c|\hat\phi(d(g))|.
$$

Since the subspace of even Paley-Wiener functions is dense in~$\CH_N^r$, we may extend the map $\phi\mapsto h_{\phi,\tau}$ by continuity provided we check that it is continuous with respect to the correct seminorms, and then the asserted formula for $\pi(h_{\phi,\tau})$ will remain valid.

By moving the contour of integration in the formula for the Fourier transform, we get $\sup|\hat\phi(x)|e^{r|x|}\le c|\phi|_{r,2}$ for some $c>0$. Thus it follows from the above estimate that for every $p>0$ there exist $r>0$, $c>0$ such that $|h_{\phi,\tau}|_{p,0,1,1}\le c|\phi|_{r,2}$.

In order to estimate the derivatives, recall from~\cite{Wa} that for every $N$ there exist $N'$ and functions $\mu$, $\nu\in C_c^N(G)$ such that
\[
h = \mu*\Delta^{N'}h + \nu*h
\]
for every $h\in\CC^1(G)$, where $\Delta$ is the Laplacian on $G$ also occurring in Proposition~\ref{convergence of spectral side}. Together with standard properties of the functions $\Xi$ and $l$ this shows that there exists $c>0$ such that
\[
|h|_{p,n,X,1}=|X\mu*\Delta^{N'}h + X\nu*h|_{p,n,1,1}\le c(|\Delta^{N'}h|_{p,n,1,1}+|h|_{p,n,1,1})
\]
for $X$ of order~$N$. We have a similar inequality for $|h|_{p,n,1,Y}$ and hence for $|h|_{p,n,X,Y}$.

Note that $\Delta$ can be chosen as $2C_K-C$, where $C_K$ is the Casimir of $K$ induced by the restriction of the Killing form of~$\g$. Specialising to $f\in\CC^p(G,\tau)$, we may replace $\Delta$ by $2d_\tau-C$, where $d_\tau=\tau(C_K)$ is a constant. Now the obvious formula
\[
(-C-b)h_{\phi,\tau}=h_{\tilde\phi,\tau},\qquad
\tilde\phi(\la)=\la^2\phi(\la)
\]
allows us to deduce the required estimate.
\qed

\section{Twisting by characters}
A {\it virtual representation} of the group $G$ is a
$\Z/2\Z$-graded finite dimensional complex representation $\psi =
\psi_+\oplus\psi_-$. We define the virtual trace and determinant
as
$$
\tr\psi(x)\= \tr\psi_+(x)-\tr\psi_-(x),\qquad \det\psi(x)\= \det\psi_+(x)/\det\psi_-(x).
$$
The function $\tr\psi(x)$ is then called a \emph{virtual character}.
Note that $\tr(\psi_1\oplus\psi_2)=\tr\psi_1+\tr\psi_2$ and
$\tr(\psi_1\otimes\psi_2)=(\tr\psi_1)(\tr\psi_2)$. Thus the virtual characters form
a ring, the \emph{character ring} ${\rm R}(G)$. 
If $G=\CG(\R)$ for a semisimple $\Q$-group as in the previous section, we call a function on $G$ parabolically regular
if it vanishes on $xP^1x^{-1}$ for every parabolic $\Q$-subgroup $\CP\ne\CG$
and every $x\in G$.

\begin{proposition}
Let $\CG$ be a connected semisimple group over $\Q$. Then there is a
virtual representation $\psi$ of $G$ such that
$\tr\psi$ is nonzero and parabolically regular.
\end{proposition}

\prf
Since characters are class functions, we may restrict attention to proper parabolic $\Q$-subgroups $\CP$ of $G$. We have the Langlands decomposition $P=MAN$, where $P^1=MN$ is defined with respect to the $\Q$-structure as in section~\ref{simple}. Since the rank of $M$ is smaller
than the rank of $G$, the restriction map ${\rm R}(G)\ra{\rm R}(M)$
is not injective. Let $\tr\psi_\CP$ be a nonzero element of its kernel. Since $N$ is
the unipotent radical of $MN$ it follows that $\tr\psi_\CP(mn)=\tr\psi_\CP(m)=0$ for all
$m\in M,\ n\in N$. Set
$$
\psi\df \bigotimes_\CP \psi_\CP,
$$
where the product runs over the finite set of all conjugacy classes of proper
parabolic $\Q$-subgroups of $G$. Then the virtual character $\tr\psi$ is
parabolically regular. Since
$\tr\psi=\prod_\CP\tr\psi_\CP$ and each $\tr\psi_\CP$ is a nonzero algebraic
function on $G$, it follows that $\tr\psi$ is nonzero.
\qed

Assume now that $\tr\psi$ is a parabolically regular virtual character.
Let $h\in C_c^\infty(G)$ and set
$$
f_\infty\= h\,\tr\psi,
$$
and let $f\fin$ be the characteristic function of~$K_\Gamma$ as in section~\ref{simple}. 

\begin{lemma}
The geometric side of the trace formula for the test function $f=f\fin\otimes f_\infty$ above is
$$
J_{\mathrm{geom}}(f) \= \sum_{[\ga]}\ \vol (\Ga_\ga \bs G_\ga)\, \tr\psi(\ga)\,\CO_\ga(h).
$$
\end{lemma}

\prf
Since $\tr\psi$ is a class function it follows that the orbital integral satisfies
$\CO_\ga(h\tr\psi)=\tr\psi(\ga)\,\CO_\ga(h)$.
\qed

This simple form of the trace formula is quite advantageous since $\tr\psi$ is
 easy to compute and concerning $h\in C_c^\infty(G)$ we have total freedom of
choice. However, the problem remains that the spectral side of the trace formula
does not simplify.
To discuss the spectral side we will consider representations $\pi\otimes\sigma$
for
$\pi\in\hat G$ and $\sigma$ finite-dimensional. We endow $\sigma$ with a $K$-invariant norm, so that these are admissible Hilbert representations. They are no
longer bounded, and hence $\pi\otimes\sigma(h)$ is not defined for $h\in L^1(G)$. However, the operator norm of $\sigma(g)$ and hence that of $\pi\otimes\sigma(g)$ grows at most exponentially with $d(g)=\dist(gK,eK)$. Thus there exists $p>0$ such that the integral defining $\pi\otimes\sigma(h)$ converges for all $h\in\CC^p_0(G)$ and $\pi\in\hat G$. Let us first prove two general facts.

\begin{lemma}\label{4.3} For a given
operator $M$ on the space of $\pi$ we
have 
$$
\tr \left(M\pi(h\tr\sigma)\right)\=\tr (M\otimes 1)(\pi\otimes\sigma(h)).
$$
\end{lemma}

\prf We compute
\eqn
\tr (M\pi(h\tr\sigma)) &=& \tr M\int_Gh(x)\pi(x)\tr\sigma(x)dx\\
&=& \tr_1\otimes\tr_2\left((M\otimes 1)\int_Gh(x)\pi(x)\otimes\sigma(x)
dx\right),
\endeqn
where $\tr_1$ and $\tr_2$ are the traces on the first and second
tensor factor.
\qed

Let $\CB(\pi)$ denote the Banach space of all bounded linear operators on the Hilbert space of $\pi$.

\begin{lemma}\label{bound norms}
For the linear map
$$
1\otimes \tr : \CB(\pi\otimes\sigma)\ra\CB(\pi)
$$
we have
$$
\|(1\otimes\tr)(T)\|\ \le\ \dim\sigma \| T \|.
$$
\end{lemma}

\prf
Let $(e_j)_j$ denote an orthonormal basis of $\sigma$. For $T\in\CB(\pi\otimes\sigma)$ and $v\in\pi$ we
have
$$
1\otimes\tr (T)(v)\=\sum_j\sp{T(v\otimes e_j)}{e_j},
$$
where the partial inner product on the right is defined as a map from $(\pi\otimes\sigma)\times
\sigma$ to $\pi$ by
$$
\sp{v\otimes w}{w'}\=\sp{w}{w'} v.
$$
This implies $\|\sp{v\otimes w}{w'}\| \le \|v\otimes w\| \|w'\|$. Hence,
\begin{eqnarray*}
\norm{1\otimes\tr(T)(v)} &\le& \sum_j\norm{\sp{T(v\otimes e_j)}{e_j}}\\
&\le& \sum_j\norm{T(v\otimes e_j)}\\
&\le & \dim\sigma\ \norm{T}\, \norm v.
\end{eqnarray*}
The lemma follows.
\qed

In order to understand the representation $\pi\otimes\sigma$ for an induced $\pi$
the following Lemma will be needed later.

\begin{lemma}\label{6.4}
Let $P=MAN$ be a parabolic subgroup of a reductive group $G$. Let $P^+=M^+AN$,
where $M^+$ is a finite index subgroup of $M$. Let
$\pi=\pi_{\xi,\nu}=\Ind_{P^+}^G(\xi\otimes\nu\otimes 1)$, where $\xi$ is an irreducible admissible representation of~$M^+$ and $\nu\in\a_\C$.
Let $\sigma$ be a finite dimensional representation of $G$. Write
$$
\sigma|_{M^+A}\=\bigoplus_{j=1}^s\sigma_j\otimes\nu_j
$$
for the decomposition into irreducibles of the restriction to the reductive group
$M^+A$.

Then, after reordering the
$\sigma_j\otimes\nu_j$ if necessary, there is a $G$-stable filtration of
$\pi\otimes\sigma$,
$$
0\= F^0(\pi\otimes\sigma)\ \subset\ \dots\ \subset\
F^s(\pi\otimes\sigma)\=\pi\otimes\sigma
$$
with quotients
$$
F^j/F^{j-1}\ \cong\ \pi_{\xi\otimes\sigma_j,\nu+\nu_j}.
$$
Since $G=P^+K$, the representation $\pi_{\xi,\nu}$ has a compact model on the
Hilbert space
$\Ind_{K\cap M^+}^K(\xi)$ which is independent of $\nu$.
In the compact model, this filtration does not depend on $\nu$.
\end{lemma}

\prf Highest weight theory implies that there is a $P^+$-stable
filtration
$$
0\= F^0\sigma\ \subset\ F^1\sigma\ \subset\ \dots\ \subset
F^s\sigma\=\sigma
$$
of $\sigma$ such that $F^j\sigma /F^{j-1}\sigma$ is isomorphic
with $\sigma_j\otimes\nu_j$ with $N$ acting trivially. Let
$\xi_\nu$ denote the representation of $P^+$ given by
$\xi_\nu(man)=a^\nu\xi(m)$.

The map $\Xi$ given by
$$
\Xi(\ph\otimes v)(x)\=\ph(x)\otimes\sigma(x)v
$$
is a $G$-isomorphism between the representation
$\pi_{\xi,\nu}\otimes\sigma$ and the induced representation
$\tilde\pi=\Ind_{P^+}^G(\xi_\nu\otimes\sigma|_{P^+})$.

Let 
$$
F^j(\tilde\pi)=\Ind_{P}^G(\xi_\nu\otimes F^j\sigma).
$$
Then the filtration $F^j(\pi\otimes\sigma)=\Xi^{-1}F^j(\tilde\pi)$
has the desired properties.
To see that the filtration does not depend on $\nu$ in the compact
model, recall that the compact model lives on the space of all
$f:K\ra V_\xi$ such
that $f(mk)=\xi(m)f(k)$ for all $m\in K\cap M^+$ and $k\in K$. Thus
$\pi\otimes\sigma$ can be modelled on the space of all $f:K\ra V_\xi\otimes
V_\sigma$ with $f(mk)=(\xi(m)\otimes 1)f(k)$. Using the construction above
it turns out that
$F^j(\pi\otimes\sigma)$ coincides with the space of all such $f$ with
$$
(1\otimes\sigma(k))f(k)\in V_\xi\otimes F^j\sigma
$$
for every $k\in K$.
\qed

\section{Convergence of the spectral side}\label{adaptMueller}

Arthur proved his trace formula for a smooth compactly
supported test function $f$ on~$\CG(\A)$. We want to substitute a
function of noncompact support depending on a
parameter, and we need uniform convergence. We are
able to prove the necessary estimates in a special case
sufficient for the purpose of this paper. Before stating
them, we should recall the definition of the spectral side
of the trace formula.

According to~\cite{Art-eisII}, Theorem~8.2, one has
$$
J_{\mathrm{spec}}(f)=\sum_\chi J_\chi(f),
$$
where $\chi$ runs through conjugacy classes of pairs
$(\CM_0, \pi_0)$ consisting of a $\Q$-rational Levi
subgroup $\CM_0$ and its cuspidal automorphic
representation~$\pi_0$, the sum being absolutely convergent.
The particular terms have expansions
$$
J_\chi(f)=\sum_{\CM,\pi} J_{\chi,\CM,\pi}(f),
$$
where the sum runs over all $\Q$-rational Levi
subgroups $\CM$ of $\CG$ containing a fixed minimal one
(which we take to be the subgroup $\CA_0$ of diagonal matrices)
and, for each~$\CM$, over all discrete automorphic representations $\pi$ of
$\CM(\A)^1$. Explicitly,
$$
J_{\chi,\CM,\pi}(f)=\sum_{s\in W_\CM}
c_{\CM,s}\int_{i(\a_{\CL}^{\CG})^*}\sum_{\CP}
\tr\left(\mathfrak M_{\CL}(\CP,\nu)M(\CP,s)
    \rho_{\chi,\pi}(\CP,\nu,f)\right)d\nu.
$$
Here, for a given element $s$ of the Weyl group of $\CM$ in $\CG$,
the Levi subgroup $\CL$ is determined by $\a_\CL=(\a_\CM)^s$,
and $\CP$ runs through all parabolic subgroups of $\CG$ having $\CM$
as a Levi component. The coefficient $c_{\CM,s}>0$ is of
no interest to us except in the case $\CM=\CG$, where it is~1.
For later use, we denote by $J_{\chi,\CM,\pi}^+(f)$ the same expression
with the trace replaced by the trace norm.

Let us comment on the items in the integrand. Let
$\rho(\CP,\nu)$ be the representation of $\CG(\A)$ which
is induced from the representation
of $\CP(\A)$ in
$$
L^2(\CM(\Q)\bs \CM(\A))\cong L^2(\CN(\A)\CP(\Q)\bs \CP(\A))
$$
twisted by~$\nu$. If one starts the induction with the
subspace of the $\pi$-isotypical component spanned by
certain residues of Eisenstein series coming from $\chi$,
one gets a subrepresentation which is denoted by
$\rho_{\chi,\pi}(\CP,\nu)$. We let $\rho_{\chi,\pi}(\CP,\nu,f)$
act in the space of $\rho(\CP,\nu)$ by composing it with the
appropriate projector. Further, there is a
meromorphic family of standard intertwining operators
$M_{\CQ|\CP}(\nu)$ between dense subspaces of
$\rho(\CP,\nu)$ and $\rho(\CQ,\nu)$ defined by an integral
for $\Re\nu$ in a certain chamber. The operator $M(\CP,s)$
is $M_{s\CP|\CP}(0)$ followed by translation with a
representative of $s$ in $\CG(\Q)$. And finally,
$\mathfrak M_{\CL}(\CP,\nu)$ is obtained from such intertwining
operators by a limiting process.

The decomposition in terms of $\chi$ is only there for technical reasons. In general, it is unknown whether the sum over $\CM$ and $\pi$ can be taken outside the sum over~$\chi$ in order to obtain an expansion in terms of the distributions
$$
J_{\CM,\pi}(f)=\sum_\chi J_{\chi,\CM,\pi}(f),
$$ 
which would be given by expressions that are analogous to $J_{\chi,\CM,\pi}(f)$ but with $\rho_{\chi,\pi}(\CP,\nu)$ replaced by
$$
\rho_\pi(\CP,\nu)=\bigoplus_\chi\rho_{\chi,\pi}(\CP,\nu).
$$
This is a problem of absolute convergence, hence of the finiteness of
\[
J_{\mathrm{spec}}^+(f)=\sum_{\chi,\CM,\pi}J_{\chi,\CM,\pi}^+(f).
\]
In~\cite{Mu2}, this problem was reduced to certain conditions on local intertwining operators, which are known to be satisfied in some cases. We will check below that those conditions are satisfied in the situation of interest to us.

Thus, we specialize to the linear algebraic
group $\CG=\SL_3$. Let $G=\CG(\R)$ be the group of real
points. We fix maximal compact subgroups $K=\mathrm{SO}_3\subset G$
and $K_p=\SL_3(\Z_p)\subset \SL_3(\Q_p)$ for all primes~$p$,
and we set $\Kfin =\prod_p K_p$. 

\begin{theorem}\label{convergence of spectral side}
Let $f=f\fin\otimes f_\infty$, where $f\fin$ is the 
characteristic function of $\Kfin $ and $f_\infty$ is a $K$-finite function in
Harish-Chandra's
$L^1$-Schwartz space $\CC^1(G)$. Then $J_{\mathrm{spec}}^+(f)<\infty$. Moreover, there exists $N>0$
with the following property. For any subset $\Pi$ of pairs $(\CM,\pi)$ as above,
there exists $c>0$ such that
$$
\sum_{(\CM,\pi)\in\Pi}J_{\CM,\pi}^+(f)\le
   c\sup_{(\CM,\pi)\in\Pi}\sup_{\nu\in i\a_\CM^*}\sup_{\CP}
   \|\rho_\pi(\CP,\nu,(1+\Delta)^N f)\|.
$$
Here we have fixed a $K$-invariant
norm on the dual of the real Lie algebra of $G$ and denoted
by $\Delta$ the corresponding element of the universal
enveloping algebra. The superscript $+$ indicates that the trace has been replaced by the trace norm.
\end{theorem}

\prf Our first assertion, which concerns absolute convergence, would follow from
Theorem~0.2 of~\cite{Mu2} if we could verify conditions 1) and $2')$ of that
theorem. Once this is done, our second assertion will be a byproduct of the proof.
Indeed, in the course of proving Lemma~6.2 of~\cite{Mu2}, equation (6.15) was used
to estimate the operator norm of $\rho_{\chi,\pi}(\CP,\nu,(1+\Delta)^Nf)$ in
terms of the $L^1$-norm of $(1+\Delta)^Nf$. If we omit that step and consider only
the terms with $(\CM,\pi)\in\Pi$, the asserted bound will follow.

The aforementioned condition~1) is a uniform bound on the derivatives of the local intertwining operators $R_{\CQ|\CP}(\pi_p,\nu)_{K_p}$ for all automorphic representations $\pi=\pi_\infty\otimes\bigotimes_p \pi_p$ of~$\CM(\A)^1$, all primes~$p$ and all open compact subgroups $K_p$ of $\CG(\Q_p)$, where the subscript $K_p$ indicates restriction to the subspace of $K_p$-fixed vectors of the representation induced from~$\pi_p$. In contrast to Theorem~0.2 of~\cite{Mu2}, our claim concerns only a fixed test function $f\fin$ which is biinvariant under a particular maximal compact subgroup~$K_p$, and hence we need only verify condition~1) for that group. However, $K_p$ as chosen above is hyperspecial, and $R_{\CQ|\CP}(\pi_p,\nu)_{K_p}$ is the identity, so that the condition is automatically satisfied.

Condition~$2')$ is a uniform bound on the derivatives of the local intertwining operators $R_{\CQ|\CP}(\pi_\infty, \nu)_\tau$ for all $\pi$ as above and all $K$-types $\tau$, where the subscript $\tau$ indicates restriction to that $K$-type in the representation induced from~$\pi_\infty$. Since we consider a fixed $K$-finite function $f_\infty$, we need only verify the condition for finitely many $K$-types, and the uniformity of the required bound in $\tau$ is no issue. However, the bound does have to be uniform in~$\pi_\infty$, which still allows us to split the set of those representations $\pi_\infty$ into a finite number of subsets and check the condition for each of them. For the subset of tempered representations, condition $2')$ follows from results of Arthur (cf.~\cite{Mu2}, Proposition~6.4).

For our group $\CG=\SL_3$, we have either $\CM=\CA_0$ or
$\CM=\CG$ or $\CM\cong \mathrm{S}(\GL_2\times \GL_1)$. In the first
case, $\pi_\infty$ is just a character of $\CA_0(\R)$,
hence tempered. In the second case $\CM=\CG$, the induced
representation $\rho_\pi(\CG)$ coincides with $\pi$, and
the intertwining operator is the identity, so that
condition $2')$ is trivially satisfied.

This leaves us with the case of the intermediate Levi
subgroups. The map $g\mapsto(g,\det g^{-1})$ is a
$\Q$-rational isomorphism from $\GL_2$ to $\mathrm{S}(\GL_2\times
\GL_1)$. Thereby we may identify $\CM$ with $\GL_2$, $K\cap
M$ with $\mathrm{O}_2$ and $K_p\cap \CM(\Q_p)$ with $\GL_2(\Z_p)$.
Thus, our only remaining concern are the automorphic
representations $\pi$ of $\CM(\A)^1\cong \GL_2(\A)^1$
occurring in $L^2(\CM(\Q)\bs \CM(\A)^1)$ for which
$\pi_\infty$ is non-tempered and which have a $\Kfin \cap
\CM(\A)$-fixed vector. For such a representation,
$\pi_\infty$ must occur in the space of $L^2$-functions on
$$
\CM(\Q)\bs \CM(\A)^1/\Kfin \cap \CM(\A) \cong \GL_2(\Z)\bs
\GL_2(\R)^1,
$$
where the last isomorphism of right $\GL_2(\R)^1$-spaces is
due to the fact that $\Q$ has class number one. The
superscript 1 refers to the subgroup of elements with
determinant of absolute value 1. Since $\GL_2(\Z)$ contains
elements with determinant~$-1$, this quotient is
isomorphic to $\SL_2(\Z)\bs \SL_2(\R)$ as an
$\SL_2(\R)$-space. As $\pi_\infty$ is non-tempered, its
Casimir eigenvalue does not exceed $1/4$ in the usual
normalization. It follows from Roelcke's eigenvalue
estimate~\cite{Ro} (cf. \cite{Hej}, ch.~11, Prop.~2.1) that
$\pi_\infty$ lies in the space of constants, so it must be
the trivial representation.

For a single representation $\pi_\infty$ and $K$-type~$\tau$,
the norm of the derivative of
$R_{\CQ|\CP}(\pi_\infty,\nu)_\tau$ is certainly bounded by a
polynomial in $\|\nu\|$ for $\nu$ outside a sufficiently
large compact subset $\Omega$ of the line $i\a_\CM^*$,
because the operator is a rational function of $\nu$.
However, the function in question is smooth and therefore
bounded on~$\Omega$, so condition~$2')$ is satisfied.
\qed

\section{Orbital integrals}

We will continue to focus on the group $G=\SL_3(\R)$.
For later use we will fix some notation. Let $P_0=M_0A_0N_0$ be the minimal
parabolic subgroup of all upper triangular matrices in $G$. We fix $A_0$ to be the
group of all diagonal matrices with positive entries and determinant one. Let
$P_1=M_1A_1N_1$ be the parabolic subgroup of all matrices in $G$ with last row of
the form $(0,0,*)$. We fix $A_1$ to be the group of all diagonal matrices of the
form $\diag(a,a,a^{-2})$, $a>0$. The group $M_1$ is isomorphic to the group
$\GL_2(\R)^1$ consisting of all real two by two matrices with determinant equal
to $\pm 1$.

For $j=0,1$ let $\rho_j\in\a_0^*$ be the modular shift of $P_j$. So by definition
for $a\in A_0$ we have $\det(a | \n_j)=a^{2\rho_j}$, where $\n_j$ is the Lie
algebra of $N_j$. Further let $\rho_{M_1}\in\a_0^*$ be the modular shift of the
parabolic $P_0\cap M_1$. Then by definition $\det(a|\n_0\cap \m_1)=a^{2\rho_{M_1}}$.
Note that $\rho=\rho_1 +\rho_{M_1}$.

The Killing form $B$ on the real Lie algebra $\g_\R=sl_3(\R)$ is
given by
$$
B(X,Y)\=\tr\ad(X)\ad(Y)\= 6\tr XY
$$
for $X,Y\in\g_\R$. We will use the same letter for its complexification as a
symmetric bilinear form on $\g$ as well as for the corresponding quadratic
form $B(X)=B(X,X)$. Let $\theta$ be the Cartan involution fixing
$K$. Then $\theta(x)= \hspace{2pt} ^t\hspace{-2pt}x^{-1}$ for $x\in G$ and $\theta(X)=- ^tX$ for $X\in\g_\R$. The Killing form gives a natural identification between $\a$ and its dual $\a^*$.
Thus it also gives an invariant form on the latter space. Note that the sum 
$\rho=\rho_1 +\rho_{M_1}$ of the last paragraph is orthogonal with respect to $B$.
Therefore, $B(\rho)=B(\rho_1) +B(\rho_{M_1})$.

We are going to apply Proposition~\ref{test} with a special choice
of the number $b$ and the representation~$\tau$ of $K\cong\mathrm{SO}_3$.
Recall that for each
$k=0,1,2,\dots$ there is an irreducible representation
$\delta_{2k}$ of dimension $2k+1$ and that this exhausts the set $\hat
K$ of irreducible representations of $K$ up to equivalence. 
For a virtual $K$-representation $\tau=\tau_+-\tau_-$ we define
$h_{\phi,\tau}=h_{\phi,\tau_+}-h_{\phi,\tau_-}$.
We choose the virtual representation
$$
\tau_0\= \delta_4-\delta_2-2\delta_0
$$
of $K$, fix $b=B(\rho_1)=1/4$ and set $h_\phi=h_{\phi,\tau_0}$. The reason for
this choice will become transparent in the next lemma.

We want to describe the action of our kernel in representations $\pi_{\xi,\nu}$ of
the principal series of~$G$. Here $P=MAN$ is a parabolic subgroup, $\xi$ a
representation in the discrete series of $M$ and $\nu\in\a^*$.
Following~\cite{MS-eta}, we write $\xi=\Ind_{M^+}^M\omega$, where $M^+$ is a
certain subgroup of finite index in~$M$ and $\omega$ belongs to the discrete series
of~$M^+$. In our situation, $M_1^+=M_1^0\cong \SL_2(\R)$ and $M_0^+=M_0$.

A $K$-finite function $f\in L^1(G)$ is called a \emph{pseudo-cusp form} if
$\tr\pi(f)=0$ for every $\pi\in\hat G$ which is induced from the minimal parabolic
$P_0$.

Let $\g=\k\oplus\p$ be the Cartan decomposition of the Lie algebra $\g$, i.e. $\p$ is the orthocomplement of $\k={\rm Lie}_\C(K)$ with respect to the Killing form.

\begin{lemma}\label{princip.series}
 Let $\pi=\pi_{\xi,\nu}$ with $\xi=\Ind_{M^+}^M\omega$ as above.
\begin{enumerate}
\item[\rm(i)] If $P=P_1$ and we denote $\p_{M_1}=\p\cap\m_1$, then
$$
\tr\pi(h_\phi)=\phi\left(|\nu|\right)
\dim\Hom_{K\cap M_1^0}\left(\omega,
{\textstyle\bigwedge^{\mathrm{odd}}\p_{M_1}-\bigwedge^{\mathrm{even}}\p_{M_1}}\right).
$$
Note in addition  that 
$$
{\textstyle\bigwedge^{\mathrm{odd}}\p_{M_1}-\bigwedge^{\mathrm{even}}\p_{M_1}}\ \cong\
(S^+-S^-)\otimes (S^--S^+)
$$
as a representation of the spin group $\Spin(B|_{\p_{M_1}})$, where $S^\pm$ are the
half-spin representations.
\item[\rm(ii)] If $P=P_0$, then $\tr\pi(h_\phi)=0$, so $h_\phi$ is a pseudo-cusp form.
\end{enumerate}
\end{lemma}

\prf Since $\pi$ can be considered as induced from~$M^+AN$, Frobenius reciprocity implies that
\begin{eqnarray*}
\tr\pi(h_\phi)&=&\phi\left(\sqrt{-\pi(C)-B(\rho_1)}\right)\dim\Hom_K(\tau_0,\pi)\\
&=&\phi\left(\sqrt{-\pi(C)-B(\rho_1)}\right)\dim\Hom_{K\cap M^+}(\omega,\tau_0)
\end{eqnarray*}
in terms of the virtual dimension of a $\Z/2\Z$-graded vector space. In case~(i), we have $K\cap M_1^0\cong \mathrm{SO}_2$, and it is straightforward to check that
$$
\tau_0|_{K\cap M_1^0}\cong{\textstyle\bigwedge^{\mathrm{odd}}\p_{M_1}-\bigwedge^{\mathrm{even}}\p_{M_1}}.
$$
If $\Lambda\in\b^*$ is the infinitesimal character of $\xi$, then $\Lambda+\nu$ is the infinitesimal character of~$\pi$, hence
\[
\pi(C)=B(\Lambda+\nu)-B(\rho)=B(\Lambda)-|\nu|^2-B(\rho_1)-B(\rho_{M_1}),
\]
as $\nu$ is imaginary and orthogonal to~$\Lambda$.
Lemma 2.4 of \cite{MS} says that $\tr\pi(h_\phi)$ vanishes unless
$B(\La)=B(\rho_{M_1})$, so we get the asserted formula.

For the proof of~(ii), observe that
$$
M_0=\{\diag(\eps_1,\eps_2,\eps_3)\mid \eps_i\in\{\pm1\},\,\eps_1\eps_2\eps_3=1\},
$$
a group whose characters are $\xi_i(\diag(\eps_1,\eps_2,\eps_3))=\eps_i$, $i=1,2,3$, together with the trivial character~$\xi_0$. It is clear that
\[
\delta_0|_{M_0}=\xi_0,\qquad
\delta_2|_{M_0}=\xi_1+\xi_2+\xi_3,\qquad
\delta_4|_{M_0}=2\xi_0+\xi_1+\xi_2+\xi_3
\]
(e.g., by dimensional reasons, using the fact that the $\xi_i$ are conjugate under the normalizer of $M_0$
in~$K$).
This implies that $\dim\Hom_{M_0}(\xi_i,\tau_0)=0$. 
\qed

Let $g\in G=\SL_3(\R)$ be regular. Then the centralizer $G_g$ is a
maximal torus of $G$, so either it is conjugate to the torus of
all diagonal elements of $G$ or to $H=A_1B$, where $B\cong \mathrm{SO}_2$
is the compact maximal torus of $M_1$. We say that $g$ is of
{\it splitrank} $2$ in the former case and of splitrank $1$ in the latter.
If $g$ is of splitrank $1$, then there is $a_gb_g\in A_1B$ and $x\in G$
such that $g=xa_gb_gx^{-1}$. Here $a_g$ is uniquely determined,
and so is $xa_gx^{-1}$, the {\it split part} of $g$.
For $g\in G$ of splitrank one with split part $\exp X$, we define its
\emph{length} by $l(g)=\sqrt{B(X)}$.
Since $\theta$ acts trivially on~$A_1$, we have $l(g)=d(a_g)$.

We will follow the conventions of \cite{HC-HA1} about the normalization of Haar-measures
on $G$ and its Lie subgroups. For $g\in G$, we use the notation
\[
D(g)=\det(\mathrm{Id}-\Ad(g^{-1})|\,\g/\g_g).
\]

\begin{proposition}\label{orbint}
Let $g\in G$ be regular, and let $\phi\in\CH_N^r$, where $r$ and $N$ are such that $h_\phi\in\CC^1_0(G)$. If the splitrank of $g$ is $2$ then
$\CO_g(h_\phi)=0$. If the splitrank of $g$ is $1$, then
$$
\CO_g(h_\phi)\= |D(g)|^{-1/2}\hat\phi(l(g)).
$$
\end{proposition}

\prf
If $\phi$ is a Paley-Wiener function, then $h_\phi\in C_c^\infty(G)$ and in
particular, $h_\phi\in \CC^2(G)$. Using Lemma \ref{princip.series} the proposition
follows from Lemma 4.3 of \cite{MS-eta}. For the general case one approximates
$\phi$ by Paley-Wiener functions.
\qed

Proposition~\ref{test} can be applied to the function
$$
\phi_\la^N(x)=(N-1)!\,(x^2+\la)^{-N},
$$ 
in which case the operator
$\phi_\la^N(\sqrt{D_\tau})$ equals $(N-1)!$ times the resolvent
$(D_\tau+\la)^{-N}$. The corresponding convolution kernel $h_{\phi,\tau}$ will be denoted by $R_{\la,\tau}^N$.
If $|\Im\sqrt{-\la}|>r$, then $\phi_\la^N\in\CH_{2N}^r$. Thus, if $|\Im\sqrt{-\la}|$ and $N$ are large enough (e.g., if $N$ and $\la>0$ are large enough), then 
 $R_{\la,\tau}^N\in\CC^1_0(G)$ and, for all for $\pi\in\hat G$,
$$
\pi(R_{\la,\tau}^N)\= (N-1)!\frac{(-\pi(C)-B(\rho_1)+\la)^{-N}}{\dim\tau}
P_{\pi,\tau}.
$$
Again, for a virtual representation $\tau=\tau_+-\tau_-$ we set $R_{\la,\tau}^N=R_{\la,\tau_+}^N-R_{\la,\tau_-}^N$, in particular $R_{\la}^N=R_{\la,\tau_0}^N$.

\begin{proposition}\label{4.1}
Let $\la>0$ and $N$ be large as above. If the splitrank of $g$ is $2$ then
$\CO_g(R_\la^N)=0$. If the splitrank of $g$ is $1$, then
$$
\CO_g(R_\la^N)\=\frac{1}{|D(g)|^{1/2}}\left(-\frac\partial{\partial\la}\right)^{N-1} 
 \frac{e^{-l(g)\sqrt{\la}}}{2\sqrt{\la}}.
$$
These orbital integrals are real and positive.
\end{proposition}

\prf
The convergence follows from the fact that $R_\la^N\in\CC_0^2(G)$ for large $\la$ and~$N$. The formula from Prop.~\ref{orbint} can be specialised using
$$
\widehat{\phi_\la^N}(x)\= \left(-\frac\partial{\partial
\la}\right)^{N-1} \frac{e^{-|x|\sqrt\la}}{2\sqrt\la}.
$$
If $\la$ is real, then $\widehat{\phi_\la^N}$ is positive, because one can see
by induction that
$$
\widehat{\phi_\la^N}(x)\=|x|^{2N-1}p_N\left(\frac1{|x|\sqrt\la}\right)e^{-|x|\sqrt\la}
$$
for some polynomial $p_N$ of degree $2N-1$ with nonnegative coefficients.
\qed

\section{Choice of the twisting character}
\label{sec3}

For the group $\SL_3(\R)$ we will now give an explicit example of a virtual
character which is parabolically regular. For this let
$\st: G\ra
\GL_3(\C)$ denote the standard representation of
$G$, and let $\eta=S^2(\st)$ be its symmetric square. Then the dimension of
$\eta$ is 6. Consider the virtual representation
$$
\psi\=\sum_{j=0}^6 (-1)^{j}\, \textstyle\bigwedge^j\eta.
$$ 
For $x\in G$ we have $\tr\psi(x)=\det(1-\eta(x))$.

\begin{lemma}\label{p-regular}
Let $P=MAN$ be a proper parabolic subgroup of $G$. Then $\tr\psi(mn)=0$ for
every $mn\in MN$. So $\tr\psi(x)$ is parabolically regular.
\end{lemma}

\prf
For $x\in G$ we have $\tr\psi(x)=0$ if $\eta(x)$ has $1$ as an eigenvalue. Since $\eta(x)$ is the
symmetric square of $x$ it has $1$ as an eigenvalue if $x$ has $\pm 1$ as an eigenvalue. So we have to
show this for $x=mn$.  There are three conjugacy classes of proper parabolic subgroups of $G$. The first is
given by the group $P_0=M_0A_0N_0$. The group $M_0N_0$ is the
subgroup of all upper triangular matrices with $\pm 1$ on the diagonal. The claim follows. The next
parabolic is $P_1=M_1 A_1 N_1$. The group $M_1 A_1$ is the centralizer of $A_1$, so
$M_1$ is isomorphic to the group of real two by two matrices of determinant $\pm 1$.
Again the claim follows. The third parabolic $P_2$ is obtained from $P_1$ by reflection along the second
diagonal. The lemma is now clear.
\qed

We now consider the Cartan subgroup $H_0$ of
diagonal matrices in~$G$. Then its Lie algebra $\a_0$, which equals the Lie algebra
of the connected component $A_0$, is the Lie algebra of diagonal matrices in~$\g$,
i.e.
$$
\a_0 = \{\diag(a,b,c) \mid a,b,c\in\C,\ a+b+c=0\}.
$$
The weight lattice is generated by
$$
\la_1(\diag(a,b,c))\= a,\ \ \ \la_2(\diag(a,b,c))\= b.
$$
The weights of the standard representation $\st$ are $\la_1,\la_2,-\la_1-\la_2$. The dominant weights are
those of the form $a\la_1+ b\la_2$ with $a,b\in\Z$, $a\ge b\ge 0$. For instance,
the modular shift of the minimal parabolic $P_0$ is
$\rho=2\la_1+\la_2$, and the modular shift $\rho_1\in\p_1^*$ of $P_1$, restricted to $\a_0$, is $\rho_1=\frac 32(\la_1+\la_2)$. The modular
shift of the minimal parabolic $P_0\cap M_1$ of $M_1$, extended trivially to~$\a_1$, is
$\rho_{M_1}=\frac 12(\la_1-\la_2)$.

For a dominant weight $\la$ let
$W_\la$ denote the irreducible representation of $G$ with highest weight $\la$. Then
$W_0$ is the trivial representation and $\st=W_{\la_1}$.

\begin{lemma}\label{psi}
The virtual representation $\psi$ decomposes as
$$
\psi\= 2W_0 - W_{2\la_1} -W_{2\la_1+2\la_2}
+W_{3\la_1+\la_2} +W_{3\la_1+2\la_2}
-W_{3\la_1}-W_{3\la_1+3\la_2}.
$$
\end{lemma}

Note that $W_{3\la_1+\la_2}$ and $W_{3\la_1+2\la_2}$ are dual to each other as are 
$W_{3\la_1}$ and $W_{3\la_1+3\la_2}$.

\prf
The representation $\eta$ contains the highest weight $2\la_1$. Weyl's dimension formula shows that
$\dim W_{2\la_1}=6=\dim\eta$, hence
$$
\eta\=W_{2\la_1}.
$$
Similarly we get $\bigwedge^2\eta = W_{3\la_1+\la_2}$ and
$\bigwedge^3\eta\=W_{3\la_1}+W_{3\la_1+3\la_2}$. Since $\bigwedge^6\eta$ is the trivial representation, the wedge product induces dualities $\eta^*=\bigwedge^5\eta$ and $\bigwedge^2\eta^*=\bigwedge^4\eta$. The automorphism $g\mapsto{}^tg^{-1}$ interchanges each of the representations $\st$ and $\bigwedge^j\eta$ with its contragredient. Its composition with a suitable Weyl group element preserves dominant weights and interchanges $\la_1$ with $\la_1+\la_2$.
\qed

\begin{lemma}\label{asymp of eta}
For elements $g\in G$ of splitrank one we have $\det(1-\eta(g))<0$ and
$$
\det(1-\eta(g))\ \sim\ -e^{l(g)} \qquad\mbox{as $l(g)\to\infty$.}
$$
\end{lemma}

\prf
For $a$, $b$, $c\in\C^\times$ we have
$$
\det(1-\eta(\diag(a,b,c)))=(1-a^2)(1-b^2)(1-c^2)(1-ab)(1-bc)(1-ca).
$$
An element $g\in G$ of splitrank one is conjugate in $\SL_2(\C)$ to the
element $\diag(e^{r+i\theta},e^{r-i\theta},e^{-2r})$, for which $l(g)=6|r|$.
The claim follows by inspection. 
\qed

\section{The geometric side}
Let $\psi$ be the virtual representation of $G$ given in
Section \ref{sec3}. Let $\phi$ be a Paley-Wiener function and define
$$
f_{\infty}^\phi(x)\= h_\phi(x) \tr\psi(x).
$$
Let $f_\phi\=f_\infty^\phi\otimes f\fin$, where $f\fin$ is the characteristic 
function of $\Kfin =\prod_p\SL_3(\Z_p)$. Let $\CE(\Ga)$ denote the set of all
conjugacy classes $[\ga]$ in $\Ga$ which are of split rank one. Note
that, according to our definition, such $\gamma$ are regular.

\begin{proposition}\label{geometric side}
The geometric side of the trace formula for $f_\phi$ is
$$
\sum_{[\ga]\in\CE(\Ga)} \vol(\Ga_\ga\bs G_\ga)
\CO_\ga(h_\phi)\det(1-\eta(\ga)).
$$

\end{proposition}

\prf It follows from Lemma~\ref{p-regular} that
$f_\infty^\phi$ is parabolically regular, so by Corollary \ref{geomtrace}
 the geometric side of the trace formula takes the
form
$$
\sum_{[\ga]}\vol(\Ga_\ga\bs G_\ga) \CO_\ga(f_\infty^\phi).
$$
The orbital integral of $f_\infty^\phi$ can be computed as
$$
\CO_\ga(f_{\infty}^\phi)\=\CO_\ga(h_\phi\tr\psi)\=\CO_\ga(h_\phi)\tr\psi(\ga),
$$
since $\tr\psi$ is invariant under conjugation. It remains to show that the sum can be reduced to the regular
classes of splitrank one. For this let $\ga\in\Ga$ with $\tr\psi(\ga)\ne 0$. Then, by the proof of Lemma~\ref{p-regular}, $\ga$ does not have $\pm1$ as
an eigenvalue.

\begin{lemma}\label{7.2}
Let $\ga\in \SL_3(\Z)$. Suppose $\ga$ does not have $\pm1$ as an eigenvalue. Then $\ga$ is regular and the $\Q$-subalgebra
$\Q(\ga)$ generated by~$\ga$ is a cubic field equal to the centralizer of $\ga$ in $\Mat_3(\Q)$. This field is complex iff $\ga$ has split rank~1.
\end{lemma}

\prf Suppose that $\gamma$ has a rational eigenvalue $\nu$.
Since the characteristic polynomial is monic and has 
integer coefficients, $\nu$ is an algebraic integer. Being rational, 
$\nu$ must be an integer, and since $\gamma^{-1}$ also has integer 
coefficients, $\nu=\pm1$.

If we exclude this case, then the characteristic polynomial of $\gamma$ is
irreducible, hence its roots are distinct and $\Q(\ga)$ is a cubic field.
Complexification shows that the centralizer $F$ of $\gamma$
in $\mathrm{Mat}_3(\mathbb Q)$ is three-dimensional and commutative.
By comparison of degree we see that $\mathbb Q(\gamma)=F$.

If $\nu$ is an eigenvalue of $\ga$, then the map which assigns
to any element of $F$ its eigenvalue in the $\nu$-eigenspace of~$\ga$ is an 
isomorphism of $F$ onto $\mathbb Q(\nu)$. Since $\ga$ has split rank one
iff it has only one real eigenvalue, the Lemma follows.
\qed

\begin{proposition}\label{convergence of trace formula}
Let $\phi=\phi_\la^N$, so that $h_\phi=R_\la^N$.
For $N,\la\gg0$ the trace formula is valid for the resulting test function~$f_\la^N$,
and Proposition~\ref{geometric side} remains valid.
\end{proposition}

\prf 
Given $N'$ and $r>0$, choose $N$ and $\la$ sufficiently large so that $\phi_\la^N\in\CH_{N'}^r$ according to Proposition~\ref{4.1}. We want to approximate $\phi_\la^N$ in this space by a sequence of Paley-Wiener functions $\phi_n$ such that $\hat\phi_n$ does not change sign and tends to $\widehat{\phi_\la^N}$ monotonely. Thus, let $\chi\in C_c^\infty(\R)$ be even, monotonely decreasing on $\R_+$ and such that $\chi(x)=1$ for $|x|\le 1$. It is easy to check that the function $\phi_n$ whose Fourier transform is $\chi(x/n)\widehat{\phi_\la^N}(x)$ does the job. Proposition~\ref{test} now implies that, given $p>0$ and $N''$, we may choose $N$ and $\la$ such that $h_{\phi_n}$ converges to $R_\la^N$ in $\CC_{N''}^p(G)$. Since $\tr\psi(g)$ grows only exponentially with~$d(g)$, we see that, if $p$ was small enough, the sequence $h_{\phi_n}\tr\psi$ converges to $R_\la^N\tr\psi$ in $\CC_{N''}^1(G)$.

It follows from Theorem~\ref{convergence of spectral side} that $f_\infty\mapsto J_{\mathrm{spec}}(f_\infty\otimes f\fin)$ is a continuous linear functional on~$\CC^1(G)$. It is clear that any continuous linear functional on $\CC^1(G)$ extends to $\CC^1_{N''}(G)$ for sufficiently large~$N''$. Thus, if we set $f_n=f_{\phi_n}$, then $J_{\mathrm{spec}}(f_n)\to J_{\mathrm{spec}}(f_\la^N)$ as $n\to\infty$. The trace formula implies that $J_{\mathrm{geom}}(f_n)\to J_{\mathrm{spec}}(f_\la^N)$ as $n\to\infty$ as well.

By Proposition~\ref{4.1} and Lemma~\ref{asymp of eta}, all the terms on the geometric side of the trace formula
for $f_n$ have the same sign, and each term tends to the corresponding term for $f_\la^N$ monotonically. Thus, we may pass to the limit $n\to\infty$ on the geometric side by monotone convergence.
\qed

An element $\ga$ of $\Ga=\SL_3(\Z)$ is called {\it primitive}
if it is not of the form $\delta^n$ for any $\delta\in\Ga$
and any natural number $n\ne1$.
For every regular $\ga\in\Ga$ there is a primitive
$\ga_0\in\Ga$ such that $\ga=\ga_0^{\mu}$ for some
$\mu\in\N$. If $\ga$ is of splitrank two then $\ga_0$ is
uniquely determined. If $\ga$ is of splitrank one then the split
part of $\ga_0$ is uniquely determined.
The normalization of Haar measures chosen \cite{HC-HA1} implies
that for each regular $\ga\in\Ga$ of splitrank $1$ we have
$\vol(\Ga_\ga\bs G_\ga)\= l(\ga_0)$.

\begin{corollary}\label{Cor7.4}
The geometric side of the trace formula
for $f_\la^N$ equals
$$
\sum_{[\ga]\in\CE(\Ga)}
l(\ga_0) \
\frac{\det(1-\eta(\ga))}{|D(\gamma)|^{1/2}} \
\left(-\frac\partial{\partial\la}\right)^{N-1}
\frac{e^{-l(\ga)\sqrt{\la}}}{2\sqrt{\la}}.
$$
\end{corollary}

\section{The spectral side}

We have proved in Proposition~\ref{convergence of trace formula} that the spectral side of the trace formula with the test function $f_\la^N$ converges for sufficiently large positive $\la$ and~$N$. Now we want to show that it extends meromorphically as a function of~$\la$ to a sufficiently large subset of the complex plane.

In the notation of section~\ref{adaptMueller}, $J_{\mathrm{spec}}(f_\la^N)$ is a sum of terms $J_{\CM,\pi}(f_\la^N)$, where $\CM$ is a Levi $\Q$-subgroup of $\CG$ and $\pi$ is a square-integrable automorphic representation of~$\CM(\A)$. The contributions with $\CM=\CG$ are the easiest ones:
$$
J_{\CG,\pi}(f_\la^N)=\tr\pi(f_\la^N)\= \tr\pi_\infty(R_\la^N\tr\psi)\prod_p \tr\pi_p(f_p)
.
$$
For any prime~$p$, we have $\pi_{triv,p}(f_p)=1$ due to our choice of~$f\fin$, while the factor at the infinite place can be computed using Lemma \ref{4.3}, giving
\[
J_{\CG,\pi}(f_\la^N)=\tr\,\pi_\infty\otimes\psi\,(R_\la^N),\qquad \pi=\pi_{triv}.
\]
We need the explicit result only for the trivial representation $\pi_{triv}$ of $\CG(\A)$. 

\begin{lemma}
$$\textstyle
\frac 1{(N-1)!}J_{\CG,\pi_{triv}}(f_\la^N)=
-2\left(\la-\frac94\right)^{-N}-2\left(\la-\frac{49}{36}\right)^{-N}
-4\left(\la-\frac14\right)^{-N}
$$
\end{lemma}

\prf
The value in question is
\[\textstyle
\frac 1{(N-1)!}\tr\psi(R_\la^N)=\sum_\sigma[\psi:\sigma]\sum_{\tau\in\hat K}
(-\sigma(C)-B(\rho_1)+\la)^{-N}[\sigma|_K:\tau]\,[\tau_0:\tau],
\]
where $\sigma$ runs through the irreducible representations of $G$ occuring in~$\psi$. The coefficients $[\psi:\sigma]$ have been calculated in Lemma~\ref{psi}, so it remains to determine $[\sigma|_K:\tau]$ and~$\sigma(C)$. All of these numbers are unchanged if we replace $\sigma$ by $\sigma^*$.

For each dominant weight $\la$ occurring we determine the decomposition of $W_\la|_K$ from its weights and compute
$$
W_\la(C)\= B(\la+\rho)-B(\rho)
$$
using $B(a\la_1+b\la_2)=\frac 19\left( a^2+b^2- ab\right)$. The result is
\begin{align*}
W_{2\la_1}|_K&=\delta_0+\delta_4, 
   & W_{2\la_1}(C) &= \frac{10}9,\\
W_{3\la_1+\la_2}|_K &= \delta_2+\delta_4+\delta_6,
   & W_{3\la_1+\la_2}(C) &= \frac{16}9,\\
W_{3\la_1}|_K &= \delta_2+\delta_6,
   & W_{3\la_1}(C) &=  2.
\end{align*}
Since $B(\rho_{M_1})=\frac1{12}$, the lemma follows from this.
\qed

Let $\Pi_\infty(\tau_0,\psi)$ be the set of all admissible irreducible
representations $\eta$ of $G$ which are subquotients of $\pi\otimes\psi$ for some {\it nontrivial} $\pi\in\hat G$ and such that $\eta$ contains a $K$-type in $\tau_0$. Let
$$
S(\tau_0,\psi)=\{\eta(C)+B(\rho_1)\mid \eta\in\Pi_\infty(\tau_0,\psi)\}.
$$
This is a closed subset of~$\C$. Note that $B(\rho_1)=\frac14$ in our normalization. Let $\Omega(\tau_0,\psi)=\C\setminus S(\tau_0,\psi)$.

\begin{proposition}\label{6.1}
The expression $J_{\mathrm{spec}}(f_\la^N)-J_{\CG,\pi_{triv}}(f_\la^N)$ extends to a holomorphic function of $\la\in\Omega(\tau_0,\psi)$.
\end{proposition}

\prf
Let us write 
$$
J_{\mathrm{spec}}(f_\la^N)-J_{\CG,\pi_{triv}}(f_\la^N)=\sum_{(\CM,\pi)\in\Pi}J_{\CM,\pi}(f_\la^N),
$$
where $\Pi$ consists of all pairs different from~$(\CG,\pi_{triv})$. We want to apply Theorem~\ref{convergence of spectral side} to show that the integral-series on the right-hand side converges normally for $\la\in\Omega(\tau_0,\psi)$ and hence represents a holomorphic function. Thus, we have to find a uniform bound on the operator norms of $\rho_\pi(\CP,\nu,f)$ for $(\pi,\CM)\ne(\CG,\pi_{triv})$, $\nu\in i\a_\CM$, parabolics $\CP$ with Levi component~$\CM$ and $f=(\Delta+1)^Nf_\la^N$, where $\la$ runs through a compact subset of $\Omega(\tau_0,\psi)$. These operators are direct sums of copies of
\[
\Ind_{\CP(\A)}^{\CG(\A)}(\pi,\nu,f)
=\Ind_P^G(\pi_\infty,\nu,f_\infty)\otimes\Ind_{\CP(\A\fin)}^{\CG(\A\fin)}(\pi\fin,\nu,f\fin),
\]
hence have the same operator norm as the latter. By our choice of $f\fin$, the second factor is the 
projection onto the subspace of $\Kfin $-fixed vectors, hence of norm one.
This leaves us with the norm of the factor at the infinite place. Thus, we
focus attention on an irreducible component of $\Ind_P^G(\pi_\infty,\nu)$,
which is a nontrivial unitary representation of~$G$. Our assertion will be a
consequence of the following result, where we use the symbol~$\pi$ in a
different sense for simplicity of notation.

\begin{lemma}
There is a uniform bound on the operator norms of
$$
\pi((\Delta+1)^N(R_{\la,\tau}^N\tr\sigma))
$$
for all nontrivial $\pi\in\hat G$, all constituents $\sigma$ of $\psi$,
all constituents $\tau$ of $\tau_0$ and $\la$ in a compact subset of $\Omega(\tau_0,\psi)$.
\end{lemma}

\prf
Recalling that $\Delta=2C_K-C$, we have
\[
\pi((\Delta+1)^N(R_{\la,\tau}^N\tr\sigma))
=\pi(R_{\la,\tau}^N\tr\sigma)\pi(2C_K-C+1)^N.
\]
Since $R_\la^N$ is $K$-finite, here the operator $C_K$ can be estimated by a
constant, so the second factor behaves like  polynomial of degree $N$ in $\pi(C)$.
Applying Lemma~\ref{bound norms} to $T=\pi\otimes\sigma(R_{\la,\tau}^N)$, we get
\begin{eqnarray*}
\norm{\pi(R_{\la,\tau}^N\tr\sigma)} &\le& \dim\sigma\,
\norm{\pi\otimes\sigma(R_{\la,\tau}^N)}\\ &=& (N-1)!\dim\sigma\,
\norm{\pi\otimes\sigma(-C-B(\rho_1) +\la)^{-N}P_\tau}.
\end{eqnarray*}

According to \cite{speh} every nontrivial $\pi\in\hat G$ is a quotient of a
representation which is parabolically induced from a unitary representation of
a proper parabolic subgroup $P=MAN$, i.e., $\pi=\pi_{\xi,\nu}$, where
$\xi\in\hat M$ and $\nu\in i\a^*$. Of course, it suffices to consider the
standard parabolics
$P_0$,
$P_1$ and $P_2$. In the case of the maximal parabolics $P_1$ and~$P_2$, if
$\xi$ itself is parabolically induced from a unitary representation of a
proper parabolic subgroup, we use induction in stages to regard $\pi$ as
induced from the minimal parabolic~$P_0$. In this way, the complementary
series corresponds to nonunitary parameters $\nu\in\a^*$ with
$\Re(\nu)=t\rho\in\a_0^*$, where $0<t\le 1/2$.

Thus, let $\pi=\pi_{\xi,\nu}$ be an induced representation.
There is a natural $K$-stable grading $Gr^j$ of $\pi\otimes\sigma$  underlying the filtration from Lemma~\ref{6.4}. The space $Gr^j$ is defined to be the set of all $f:V\ra
V_\xi\otimes V_\sigma$ such that
$$
(1\otimes\sigma(k))f(k)\ \in\ V_\xi\otimes(\sigma_j\otimes\nu_j).
$$
Let $P_j$ be the projection to $Gr^j$. Then
$$
P_jf(k)\= (1\otimes\sigma(k^{-1}))(1\otimes
Pr_{\sigma_j\otimes\nu_j})(1\otimes\sigma(k))f(k).
$$
With respect to this grading the operator
$$
\pi\otimes\sigma(-C-B(\rho_1)+\la)^NP_\tau
$$
is a triangular matrix whose
entries are polynomials in $\nu$ of degree $\le 2N$.
The diagonal entries have the form
$$
(-\eta(C)-B(\rho_1)+\la)^N
$$
with $\eta$ being a subquotient of~$\pi\otimes\sigma$, and their leading term in $\nu$ is $B(\nu)^N$. Hence the inverse matrix is triangular and its
entries are rational functions in $\nu$ which tend to zero as fast as
$B(\nu)^{-N}$ as $\nu\ra\infty$.
Moreover, these functions have no poles at points $\nu$ parametrising $\eta\in\Pi_\infty(\tau_0,\psi)$ if $\lambda\in\Omega(\tau_0,\psi)$, as follows from the definition of the latter set.
This implies that the norm $\norm{\pi\otimes\sigma(R_{\la,\tau}^N)}$ times $(1+|\pi(C)|)^N$ is bounded.

The bound will depend on $\xi$, but for $P$ being the minimal parabolic, the group
$M$ is finite and there are only finitely many $\xi$. For the maximal
parabolics we may assume that $\xi$ is not induced itself, so it is one-dimensional or a (limit of) discrete series representation. As there are only
finitely many such $\xi$ for which some $\xi\otimes\sigma_j$ has a $K\cap M$-type in $\tau|_{K\cap M}$, we get a uniform bound on the operator norm in question, and the lemma follows.
\qed

As noted above, Proposition~\ref{6.1} is thereby proved, too. Next we show that the set $\Omega(\tau_0,\psi)$, to which we have meromorphically continued the spectral side, is large enough for our goals, in particular, that it contains the main pole at~$\la=\frac94$.

\begin{proposition}\label{6.5}
There exists $\al<\frac32$ such that $\Re\bigl( \pm \sqrt\la\,\bigr)\le\al$
for every $\la\in S(\tau_0,\psi)$.
\end{proposition}

For the proof we will need two prerequisites. First, let $S_\al$ be the set of all $z\in\C$ with $|\Re(z)|\le\al$, and let
$S_\al^2=\{ z^2\mid z\in S_\al\}$. A computation shows that
$$
S_\al^2\=\left\{ x+iy\left|\, x+\frac{y^2}{4\al^2}\le\al^2\right.\right\}.
$$

\begin{lemma}\label{7.7}
Let $V_\R$ be a finite dimensional real vector space with a positive definite
symmetric bilinear form $B_\R$. Let $V$, $B$ be their complexifications. Then every
$v\in V$ can be written as $v=\Re v +i\Im v$ for $\Re v,\Im v\in V_\R$. For every
$v\in V$ and all $\beta >0$, $c\ge 0$ we have
$$
B(\Re v)-c\in S_\beta^2\ \Rightarrow\ B(v)-c\in S_\al^2,
$$
where $\al=\sqrt{\beta^2+c}$.
\end{lemma}

\prf
We first show the case $c=0$, so $\al=\beta$.
For $v\in V$ write $B(v)=x+iy$. Then $x=B(\Re v)- B(\Im v) $ and $y=2B(\Re v,\Im
v)$. The Cauchy-Schwartz inequality implies  $y^2\le 4 B(\Re v)B(\Im v)$ and thus
$$
x+\frac {y^2}{4\al^2}\ \le\ B(\Re v)+B(\Im v)\left( \frac{B(\Re
v)}{\al^2}-1\right).
$$
Now $B(\Re v)\in S_\al^2$ implies $B(\Re v)\le\al^2$ and hence $x+\frac
{y^2}{4\al^2}\le B(\Re v)\le\al^2$, whence the claim.

For the general case simply observe that $S_\beta^2+c\subset S_\al^2\subset S_\al^2+c$.
\qed

We will also need some elementary facts about representations of $\SL_2(\R)$. For $n\in\Z$, consider the character of the subgroup of upper triangular matrices which takes value $\mathop{\rm sgn}(a)a^n$ on a matrix with upper left entry~$a$, and let $\pi_n$ be the normalised induced representation of~$\SL_2(\R)$. Then $\pi_n$ has a unique subrepresentation $\xi_n^+$ (resp.\ $\xi_n^-$) whose $\SO_2$-types are bounded only from below (resp.\ only from above). In fact,
\[
\xi_n^\pm|_{\SO_2}\cong
\bigoplus_{m=0}^\infty\eps_{\pm(n+2m+1)},
\]
where $\eps_r$, $r\in\Z$, are the characters of $\SO_2$. For $n>0$ (resp.\ $n=0$), the representations $\xi_n^\pm$ are irreducible and constitute the discrete series (resp.\ its limits). In the notation of~\cite{knapp}, $\xi_n^\pm=\CD_{n+1}^\pm$.

\begin{lemma}\label{SL2}
If $\zeta$ is a $k+1$-dimensional irreducible representation of~$\SL_2(\R)$, then $\xi_n^\pm\otimes\zeta$ has a filtration whose subquotients are isomorphic to $\xi_{n+m}^\pm$ with $|m|\le k$ and $m\equiv k\pmod2$.
\end{lemma}

\prf
By Lemma~\ref{6.4}, $\pi_n\otimes\zeta$ has a filtration with subquotients isomorphic to $\pi_{n+m}$ with $|m|\le k$ and $m\equiv k\pmod2$. This induces a filtration on $\xi_n^+\otimes\zeta$ of length at most $k+1$ whose subquotients are subrepresentations of the subquotients of the previous filtration. Since
\[
\zeta|_{\SO_2}\cong\eps_{-k}\oplus\eps_{-k+2}\oplus\dots\oplus\eps_k,
\]
the $\SO_2$-types $\eps_r$ of $\xi_n^+\otimes\zeta$ are bounded from below and have multiplicity $k+1$ for large~$r$. This identifies those subrepresentations of $\pi_{n+m}$ uniquely. The case of $\xi_n^-\otimes\zeta$ is analogous.
\qed

\noindent{\bf Proof of Proposition~\ref{6.5}:} We are going to show that $\al=\sqrt{19/12}$ has the required property, i.e., that $\eta(C)+B(\rho_1)\in S_\al^2$ for all $\eta\in\Pi_\infty(\tau_0,\psi)$.
Note that for this $\al$ we have $\al=\sqrt{\beta^2+c}$ for $\beta=\sqrt{3/2}$ and $c=1/12$.
Like in the proof of Proposition~\ref{6.1}, we will use the classification
of $\hat G$ from~\cite{speh}. Consider the representation
$\pi=\pi_{\xi,\nu}$ induced from the parabolic subgroup $P=MAN$ of~$G$,
where $\xi\in\hat M$ and $\nu\in\a^*$.
If $\eta$ is a subquotient of $\pi\otimes\sigma$, then by Lemma~\ref{6.4} there is a $j$ such that $\eta$ is a subquotient of $\pi_{\xi\otimes\sigma_j,\nu+\nu_j}$.

Let us start with the case that $\pi$ is induced from the minimal
parabolic $P_0=M_0A_0N_0$ of all upper triangular matrices,
where $M_0$ is the group of all diagonal matrices with entries $\pm 1$
and $N_0$ is the group of all upper triangular matrices with ones on the
diagonal. We have
$$
\eta(C)+B(\rho_1)\= B(\nu+\nu_j)-B(\rho_{M_1})\= B(\nu+\nu_j)-\frac 1{12}.
$$
In this case $\nu_j$ runs through the weights of~$\sigma$. For the principal series, we have $\nu\in\a_0^*$ purely imaginary. Since we discuss the complementary series as induced from~$P_0$, we also have
to consider $\nu\in\a_0^*$ with $\Re(\nu)=t\rho$, $0<t\le\frac12$. In the same
manner, we can at once handle the representations induced from the
one-dimensional representations of
$P_1$ or~$P_2$ by embedding them into principal series for the
parameter~$\frac12\rho$.

Refering to Lemma \ref{7.7}, we see that it suffices to show $B(t\rho+\nu_j)-\frac 1{12}\in S_\beta^2$ for $0\le t\le\frac12$. On a Weyl orbit of weights $\nu_j\in\a_0^*$, this function obtains its largest value for dominant weights. Those occurring in $\psi$ are of
the form $\nu_j=a\la_1+b\la_2$ for $0\le b\le a\le 3$. The maximum is obtained for $a=b=3$ and $t=\frac12$, where we get $B(3\la_1+3\la_2+\rho/2)-\frac
1{12}=\frac 32=\beta^2$.

It remains to consider the case when $\pi=\pi_{\xi,\nu}$ is induced
from a maximal parabolic, $\xi$ belongs to the (limit of) discrete series and $\nu$ is purely imaginary. Since all maximal parabolics are conjugate under the automorphism group of~$G$, it suffices to consider~$P_1=M_1A_1N_1$. Then $M_1\cong\GL_2(\R)^1$ has two connected components, and $\xi=\Ind_{M_1^0}^{M_1}(\xi^+)$, where $\xi^+$ is in the (limit of) discrete series of $M_1^0\cong \SL_2(\R)$. Since conjugation by the nontrivial element of $M_1/M_1^0$ switches holomorphic with antiholomorphic discrete series, we may assume that $\xi^+\cong\xi_n^+$ for some~$n\ge0$ in the notation of Lemma~\ref{SL2}. By induction in stages, we consider $\pi$ as being induced from $P_1^0=M_1^0A_1N_1$, so $\eta$ is a subquotient of $\pi_{\xi_n^+\otimes\zeta,\nu+\nu_j}$ for some irreducible constituent $\zeta$ of~$\sigma_j$, now denoting a representation of~$M_1^0$.

The highest weight of $\zeta$ is of the form $k\rho_{M_1}\in\a_{M_1}^*$, where $k$ is a nonnegative integer because $\rho_{M_1}=\frac12(\la_1-\la_2)$ happens to be the generator of the weight lattice of~$A_{M_1}$. Now $\dim\zeta=k+1$, and from Lemma~\ref{SL2} we conclude that $\eta$ is a subquotient of $\pi_{\xi_{n+m}^+,\nu+\nu_j}$ for some $|m|\le k$. The condition $\eta\in\Pi_\infty(\tau_0,\psi)$ means that $\eta$ has a $K$-type in common with $\tau_0$, and by Frobenius reciprocity this implies that $\xi_{n+m}^+$ has a $K_{M_1^0}$-type in common with $\tau_0|_{K_{M_1^0}}$. Since the only constituents of $\tau_0$ are $\delta_{2k}$ with $|k|\le2$ and
\[
\delta_{2k}|_{K_{M_1^0}}\ \cong\
\eps_{-k}\oplus\eps_{-k+1}\oplus\dots\oplus\eps_{k},
\]
this imposes the restriction $n+m\le1$. Remembering that $n\ge0$, we see that the infinitesimal character of $\xi_{n+m}^+$, which is $(n+m)\rho_{M_1}$, lies in the segment connecting $\rho_{M_1}$ with the lowest weight $\omega_j$ of~$\sigma_j$,
hence is of the form $u\omega_j+t\rho_{M_1}$ with $|u|\le1$ and $0\le t\le1$. 
But $\pm\omega_j+\nu_j$ is an $\a_0$-weight occurring in~$\psi$, and so the infinitesimal character of $\eta$, which is $(n+m)\rho_{M_1}+\nu+\nu_j$, can be written as $\omega+\nu+t\rho_{M_1}$, where $\omega$ is in the convex hull of the weights occurring in $\psi$ and $0\le t\le1$. Thus
\[
\eta(C)+B(\rho_1)=B(\omega+\nu+t\rho_{M_1})-B(\rho_{M_1}).
\]
If we replace $\rho_{M_1}$ by its Weyl-conjugate $\frac12\rho$, it follows from the calculation done in the case of~$P_0$ that $\eta(C)+B(\rho_1)\le\frac32=\beta^2$.
\qed

\section{The prime geodesic theorem}

Note that if $\ga\in\Ga$ is of splitrank one, then $\ga$ is not 
conjugate to its inverse $\ga^{-1}$.
So let $\CE^\pm(\Ga)$ be the set $\CE(\Ga)$ of conjugacy classes of
 split rank one modulo the equivalence relation $[\ga]\sim [\ga^{-1}]$.
Let $\CE_0^\pm(\Ga)$ be the subset of primitive classes.

For $[\ga]\in\CE_0(\Ga)$ let $N(\ga)=e^{l(\ga)}$ and define for $x>0$,
$$
\pi(x)\= \# \{ [\ga]\in\CE_0^\pm(\Ga)\mid N(\ga)\le x\}.
$$

\begin{theorem}[Prime Geodesic Theorem]\label{PGT}{\ }\\
For $x\ra\infty$ we have the asymptotic formula
$$
\pi(x)\ \sim\ \frac x{\log x}.
$$
\end{theorem}

Our main result Theorem \ref{main} can be deduced from Theorem~\ref{PGT} as follows. If
$[\ga]\in\CE_0(\Ga)$, then by Lemma 7.2 the centralizer $F_\ga$ of $\ga$ in $\Mat_3(\Q)$
is a complex cubic field and the set $\CO_\ga=F_\ga\cap\Mat_3(\Z)$ is an order in
$F_\ga$ whose unit group is generated by $\ga$. We claim that every order $\CO$
occurs $h(\CO)$ times in this way. The corresponding claim for a division algebra
instead of $\Mat_3(\Q)$ is shown in \cite{class}, section 2.
In that section the restriction to a division algebra was made to secure that the
centralizer $F_\ga$ would be a field. In the present paper this information is
obtained from Lemma \ref{7.2}. Thus the proof goes through. 
Since $\ga$ is
primitive, it is a generator of the unit group
$\CO_\ga^\times=\pm \ga^\Z$. Comparing the metric given by the Killing form
and the measure which defines the regulator \cite{neukirch}, one finds that
$l(\ga)=3R(\CO_\ga)$. Thus Theorem \ref{main} follows from the Prime
Geodesic Theorem.

It remains to prove Theorem \ref{PGT}. For $\ga\in\CE(\Ga)$ let
$$
c_\ga\= \frac{-\det(1-\eta(\ga))}{|D(\ga)|^{1/2}}e^{-l(\ga)/2}.
$$
Since $|\rho_1|=\frac12$, we have $|D(\ga)|\sim e^{l(\ga)}$ and hence
$c_\ga\ra 1$ as $l(\ga)\ra\infty$ by Lemma \ref{asymp of eta}.

\begin{proposition}\label{7.3}
Let
$$
L(s)\=\sum_{[\ga]\in\CE^\pm(\Ga)} l(\ga_0)\, c_\ga\, e^{-l(\ga)s}.
$$
where $\ga_0$ is a primitive element whose power is~$\gamma$. 
Then $L(s)$ converges for $\Re(s)>1$ and extends to a meromorphic function on
$\{\Re(s)> 1-\eps\}$ for some $\eps >0$. It has a simple pole of residue $1$ at
$s=1$ and is otherwise analytic on $\{\Re(s)\ge 1\}$.
\end{proposition}

\prf 
For $\la\gg0$ let $M_N(\la)$ be $-\frac12$ times the geometric side of
the trace formula for $f_\la^N$.
By Corollary \ref{Cor7.4},
$$
M_N(\la)\=-\sum_{[\ga]\in\CE^\pm(\Ga)} l(\ga_0)c_\ga\left(-
\frac\partial{\partial\la}\right)^{N-1}\frac{e^{-l(\ga)(\sqrt{\la}-\frac 12)
}}{2\sqrt{\la}}.
$$
For $\Re(s)\gg0$ we compute formally at first,
\begin{eqnarray*}
M_N(s(s+1)+\frac 14) &=& (-1)^N\sum_{[\ga]}l(\ga_0)c_\ga
\left(\frac\partial{\partial(s(s+1))}\right)^{N-1}\frac{e^{-sl(\ga)}}{2s+1}\\
&=& (-1)^N\sum_{[\ga]}l(\ga_0)c_\ga
\left(\frac
1{2s+1}\frac\partial{\partial(s)}\right)^{N-1}\frac{e^{-sl(\ga)}}{2s+1}\\
&=&
(-1)^N\frac 1{2s+1}D^{N-1} L(s),
\end{eqnarray*}
where $D$ is the differential operator $D=\frac\partial{\partial s}\frac 1{2s+1}$.

We deduce from the proof of Proposition~5.3 that there exists $C>0$ such 
that, for $\la>0$,
$$
\left(-\frac{\partial}{\partial\la}\right)^{N-1} 
                          \frac{e^{-|x|\sqrt\la}}{2\sqrt\la}
\ \ge\ C \frac{e^{-|x|\sqrt\la}}{\la^{N-1/2}},
$$
which shows that $M_1(\la)$ is 
convergent for $\la\gg0$.
From this and Propositions \ref{6.1} and \ref{6.5} we infer the claim on
the analytic continuation. Since $L(s)$ is a Dirichlet series with positive
coefficients by Lemma~\ref{asymp of eta} and Proposition~\ref{4.1}, its analytic continuation
also implies convergence. Proposition~\ref{7.3} is proven.
The Prime Geodesic Theorem follows from the Proposition by the Wiener-Ikehara Theorem as in \cite{chand}.
\qed

\end{document}